\documentclass[11pt]{amsart}
\usepackage{a4,amsmath,amssymb,verbatim}

\newtheorem{prop}[equation]{Proposition}

\newtheorem{thm}[equation]{Theorem}
\newtheorem{cor}[equation]{Corollary}
\newtheorem{lem}[equation]{Lemma}

\theoremstyle{definition}

\newtheorem{exa}[equation]{Example}

\numberwithin{equation}{section}

\newcommand{\sands}{\mbox{$\quad\text{and}\quad$}}

\newcommand{\Hom}{\operatorname{Hom}}

\newcommand{\CPo}{\mbox{$\text{\it CP\/}^1$}}

\newcommand{\CPn}{\mbox{$\text{\it CP\/}^n$}}
\newcommand{\CP}[1]{\mbox{$\text{\it CP\/}^{#1}$}}

\newcommand{\cs}{\mathbin{\#}}
\newcommand{\oz}{\mbox{$\overline{\zeta}$}}

\newcommand{\osu}{\mbox{$\varOmega_*^{U}$}}
\newcommand{\ous}{\mbox{$\varOmega^*_{U}$}}

\newcommand{\bZ}{\mathbb{Z}}
\newcommand{\bR}{\mathbb{R}}
\newcommand{\bC}{\mathbb{C}}

\newcommand{\bA}{\mathbb{A}}

\newcommand{\mtn}{\mbox{$M^{2n}$}}
\newcommand{\dfmtn}{\mbox{$M_\bullet^{2n}$}}

\newcommand{\brp}{\mbox{$\bR_{\scriptscriptstyle\geqslant}$}}

\newcommand{\brpf}{\mbox{$\bR_{\scriptscriptstyle\geqslant}^\mathcal{F}$}}
\newcommand{\brsf}{\mbox{$\bR_{\scriptscriptstyle >}^\mathcal{F}$}}
\newcommand{\brpg}{\mbox{$\bR_{\scriptscriptstyle\geqslant}^\mathcal{G}$}}
\newcommand{\brsg}{\mbox{$\bR_{\scriptscriptstyle >}^\mathcal{G}$}}
\newcommand{\brpn}{\mbox{$\bR_{\scriptscriptstyle\geqslant}^n$}}
\newcommand{\brsn}{\mbox{$\bR_{\scriptscriptstyle >}^n$}}
\newcommand{\wrp}{\mbox{$W^\bR(P)$}}

\newcommand{\xft}{\mbox{$X(F)^{2(n-1)}$}}

\newcommand{\xgt}{\mbox{$X(G)^{2(n-1)}$}}
\newcommand{\xgk}{\mbox{$X(G)^{2(n-k)}$}}

\newcommand{\dfxgk}{\mbox{$X(G)_\bullet^{2(n-k)}$}}

\newcommand{\nsplx}{\mbox{$\varDelta^n$}}
\newcommand{\splx}[1]{\mbox{$\varDelta^{#1}$}}

\begin{document}
\bibliographystyle{plain}
\title[Tangential Structures, Toric Manifolds \& Polytopes]
{Tangential Structures on Toric Manifolds, and Connected Sums of
Polytopes}
\author{Victor M Buchstaber}
\address{Department of Mathematics and Mechanics, Moscow State
University, 119899 Moscow, Russia}
\email{buchstab@nw.math.msu.su}
\author{Nigel Ray}
\address{Department of Mathematics, University of Manchester,
Manchester M13~9PL, England}
\email{nige@ma.man.ac.uk}


\date{13.9.00}

\begin{abstract}

We extend work of Davis and Januszkiewicz by considering {\it
omnioriented\/} toric manifolds, whose canonical codimension-2
submanifolds are independently oriented. We show that each
omniorientation induces a canonical stably complex structure,
which is respected by the torus action and so defines an element
of an equivariant cobordism ring. As an application, we compute
the complex bordism groups and cobordism ring of an arbitrary
omnioriented toric manifold. We consider a family of examples
$B_{i,j}$, which are toric manifolds over products of simplices,
and verify that their natural stably complex structure is induced
by an omniorientation. Studying connected sums of products of the
$B_{i,j}$ allows us to deduce that every complex cobordism class
of dimension $>2$ contains a toric manifold, necessarily
connected, and so provides a positive answer to the toric analogue
of Hirzebruch's famous question for algebraic varieties. In
previous work, we dealt only with disjoint unions, and ignored the
relationship between the stably complex structure and the action
of the torus. In passing, we introduce a notion of connected sum
$\#$ for simple $n$-dimensional polytopes; when $P^n$ is a product
of simplices, we describe $P^n\# Q^n$ by applying an appropriate
sequence of {\it pruning operators}, or hyperplane cuts, to $Q^n$.

\end{abstract}

\maketitle

\section{Introduction}\label{intro}

The study of toric varieties (or torus embeddings, as they were
originally known) has entranced algebraic geometers since the 1970s,
and provides a host of elegant and illuminating examples. Several
comprehensive textbooks are now available, by authors such as Ewald
\cite{ew:cca}, Fulton \cite{fu:itv}, and Oda \cite{od:cba}. In their
pioneering paper of 1991 \cite{daja:cpc}, Davis and Januszkiewicz
defined the related notion of toric manifold, thereby extending the
audience for the fascinating interplay between combinatorics, geometry
and topology which characterises the subject.

A toric manifold $\mtn$ admits a smooth action of the torus $T^n$
which may be identified locally with the standard action of $T^n$ on
$\bC^n$; the quotient space is required to be an $n$-dimensional ball,
invested with the combinatorial structure of a simple convex polytope
by the fixed point sets of appropriate subtori. A classic example is
provided by the complex projective space $\CPn$, whose quotient
polytope is the $n$-simplex $\nsplx$. More general examples may fail
to be complex, as shown by the connected sum $\CPn\cs\CPn$ whose
quotient polytope is the product of simplices
$\splx{1}\times\splx{n-1}$. However, it is clear from the work of
Davis and Januskiewicz that the action of the torus gives rise to a
family of complex structures on the {\it stable\/} tangent bundle; our
basic aim is to classify the members of this family in terms of
omniorientations, and to develop the consequences for complex
cobordism theory.

Our programme has its origins in \cite{bura:tmc}, where we
constructed a sequence of stably complex toric manifolds
$B_{i,j}$. Although we ignored any relationship between their
stably complex structure and the action of the torus, we did
confirm that the resulting cobordism classes are multiplicative
generators of the complex cobordism ring $\osu$. It follows that
every cobordism class is represented by a disjoint union of toric
manifolds, which are suitably oriented products of the $B_{i,j}$.
This situation compares with Hirzebruch's corresponding result
\cite{hi:km} for algebraic varieties, and leads to the same
question; can the representatives be chosen to be connected? For
varieties it remains unanswered, but we prove:

\smallskip

\noindent
{\bf Theorem 6.11.\;}
{\it
In dimensions $>2$, every complex cobordism class contains a toric
manifold, necessarily connected, whose stably complex structure is
induced by an omniorientation, and is therefore compatible with the
action of the torus.}

\smallskip

During our proof we develop the notion of connected sum $\#$ for simple
polytopes. In the case of the examples $B_{i,j}$ the quotient polytopes
are products of simplices; any such $P^n$ determines a sequence of {\it
pruning operators}, whose application to an arbitrary simple polytope
$Q^n$ provides an alternative description for $P^n\cs Q^n$ in terms of
hyperplane cuts.

Davis and Januszkiewicz succeeded in computing the integral homology
and cohomology of an arbitrary toric manifold $\mtn$, and our analysis
of stably complex structures allows us to extend their computations to
complex cobordism theory. We summarise this development in a separate
section, having studied the specific case of bounded flag manifolds in
\cite{bura:fml}. Our statement of results for $\mtn$ depends on the
choice of stably complex structure, and needs care to make
precise. Nevertheless, knowledge of $\ous(\mtn)$ leads to the
description of $E^*(\mtn)$ for any complex oriented cohomology theory,
and has already been translated by Strickland \cite{st:fge} into the
more sophisticated language of formal groups. Non-oriented theories
are are more difficult to deal with, but Bahri and Bendersky have made
considerable progress in the case of $KO$-theory \cite{babe:ktt}.

Recently, Hattori \cite{ha:act} and Masuda \cite{ma:utm} have studied
torus actions on stably complex manifolds as a generalisation of toric
varieties, and our work confirms that omnioriented toric manifolds fall
within their framework. Also, Panov \cite{pa:cfc}, \cite{pa:hgm} has
incorporated our notion of omniorientation into his combinatorial
description of certain cobordism invariants. Although the methods and
objectives of these authors differ significantly from ours, it will be
of interest to combine the various approaches in future.

Throughout our work we write $T^n$ for the $n$-dimensional torus,
and refer to its representation by diagonal matrices in $U(n)$ as
the {\it standard\/} action on $\bC^n$. The quotient space of this
action is the positive cone
\[
\{(x_1,\dots,x_n): x_r\geq 0,\;\text{for}\;1\leq r\leq n\}
\]
in $\bR^n$, which we write as \brpn. We may recover $\bC^n$ from the
cone as the identification space $(T^n\times\brpn)/\!\!\approx$,
where $(t,x)\approx(u,x)$ whenever the coordinates of $t$ differ
from the coordinates of $u$ only in those positions $r$ where $x_r$
is zero, and the standard action is induced by multiplication in
$T^n$. We let $\brsn$ denote the subspace of vectors whose
coordinates are strictly positive. On several occasions we consider
smooth manifolds which are locally diffeomorphic to \brpn; following
the original definitions of the 60s \cite{ja:com}, we refer to these
as $n$-dimensional manifolds {\it with corners}.

We often abbreviate singleton sets such as $\{v\}$ by omitting the
brackets.

It is a pleasure to acknowledge the insight we have gained from
discussions with several of our colleagues during the preparation of
this work, and to nominate Tony Bahri, Yusuf Civan, Taras Panov and
Neil Strickland for particular mention. The second author also wishes
to apologise for the long delay in completing the manuscript; this was
beyond his control, but has led to some confusion over the status of
our results.

\section{Toric manifolds}\label{toma}

We begin with a summary of Davis and Januszkiewicz's treatment of
toric manifolds, giving our own interpretation as required for the
study of stably complex structures in later sections.

We consider an unordered set $\mathcal{H}$ of $m$ closed halfspaces $H$
in a real affine space $\bA^n$. We assume that $m>n$, and that the
bounding hyperplanes are in general position, so that no $n+1$ of them
meet; we also insist that the removal of any halfspace will enlarge the
intersection $\cap_\mathcal{H}H$. We refer to $\cap_\mathcal{H}H$ as a
{\it simple $n$-polyhedron} $P^n$, and to each of its intersections with
a bounding hyperplane as a {\it facet}. Once an orthonormal coordinate
system is chosen for $\bA^n$ we may represent $P^n$ by a matrix
inequality $A_Px\geq b$, where $A_P$ is an $\mathcal{H}\times n$ real
matrix, and $x$ and $b$ are column vectors in $\bR^n$ and
$\bR^\mathcal{H}$ respectively; the rows of $A_P$ are indexed by the
elements of $\mathcal{H}$, so that each $H$ is described by the
corresponding row of the inequality. We reserve the term {\it
polytope\/} for a bounded polyhedron.

The standard octahedron, for example, is not simple, because $4$-tuples
of bounding hyperplanes meet at four of its vertices in $\bR^3$.

A simple polyhedron $P^n$ is uniquely determined by its set of facets
$\mathcal{F}(P)$, which we abbreviate to $\mathcal{F}$ whenever
possible. Given $0<k\leq n$, every nonempty intersection of $k$ facets
forms a {\it face\/} of $P^n$, which has codimension $k$ by general
position; conversely, every codimension-$k$ face $G$ determines a
unique set $\mathcal{F}_G$ of $k$ facets. Any such face is itself a
simple polyhedron $G^{n-k}$, defined by those facets of $P^n$ which
intersect it properly. We may therefore partition $\mathcal{F}$ as
\begin{equation}\label{partfacs}
\mathcal{F}(G)\cup\mathcal{F}_G\cup\mathcal{D}_G,
\end{equation}
where $\mathcal{D}_G$ consists of the facets disjoint from $G$. In
particular, every vertex $v$ of $P^n$ is determined by a unique set
$\mathcal{F}_v$ of $n$ facets, and so lies in a neighbourhood which is
linearly isomorphic to the cone \brpn. It follows that $P^n$ is an
$n$-dimensional manifold with corners, and has an atlas with one affine
chart $U_v$ for each vertex $v$. Clearly $P^n$ is a convex submanifold
of the ambient $\bA^n$.

Geometers originally studied polyhedra up to affine equivalence, but the
weaker notion of {\it combinatorial equivalence}, determined by the
lattice of faces $\mathfrak{L}_F(P)$, is now equally fashionable;
polyhedra are combinatorially equivalent if and only if they are
diffeomorphic as manifolds with corners. Many fascinating details, and a
host of further references, are given in Ziegler's book \cite{zi:lp}.

To establish our notation we describe two fundamental families of
polytopes in some detail, namely simplices and cubes.

The standard $n$-simplex $\splx{n}$ lies in $\bR^n$, and has defining
halfspaces 
\begin{equation}\label{nsplx}
H_r=\{x:x_r\geq 0\}\;\;\text{for $1\leq r\leq n$, and}\;\;
H_{n+1}=\{x:x_1+\dots+x_n\leq 1\},
\end{equation}
with corresponding facets $D_r=\splx{n}\cap H_r$. Each $D_r$ is a
copy of the $(n-1)$-simplex $\splx{n-1}$ for $1\leq r\leq n$, whilst
$D_{n+1}$ is affinely equivalent to $\splx{n-1}$. The codimension
$k$ faces $D_S$ are $(n-k)$-simplices, indexed by the $k$-element
subsets $S$ of $\{1,\dots,n+1\}$, and the face lattice
$\mathfrak{L}_F(\splx{n})$ is therefore Boolean of rank $n$.

The standard $n$-cube $I^n$ also lies in $\bR^n$, and has defining
halfspaces 
\begin{equation}\label{ncube}
H_r^0=\{x:x_r\geq 0\}\sands H_r^1=\{x:x_r\leq 1\}
\end{equation}
for $1\leq r\leq n$, with corresponding facets
$C_r^\varepsilon=I^n\cap H_r^\varepsilon$, where $\varepsilon=0$ or
$1$. Each $C_r^\varepsilon$ is an $(n-1)$-cube $I^{n-1}$, for $1\leq
r\leq n$. The codimension $k$ faces are $(n-k)$-cubes, indexed by
the cartesian coordinates of their centres; these are ternary
sequences $\xi$ of length $n$ on $\{\frac{1}{2},0,1\}$, in which
$\frac{1}{2}$ occurs $n-k$ times. Thus $C_r^\epsilon$ is indexed by
$\xi_j=\frac{1}{2}$ for $j\neq r$ and $\xi_r=\epsilon$, whilst the
vertices are given by their coordinate sequences of $0$s and $1$s.
The face lattice $\mathfrak{L}_F(I^n)$ has $3^n$ elements, and is of
independent interest to combinatorialists.

We often use the product polytope $I^m\times\nsplx$. This has facets
$C_r^\varepsilon\times\nsplx$ and $I^m\times D_s$ for $1\leq r\leq
m$ and $1\leq s\leq n$, written $E_r^\varepsilon$ and $E_s$
respectively.

We consider $2n$-dimensional manifolds $\mtn$ which are equipped with an
action $\alpha$ of the torus $T^n$, and suppose for convenience that
both $\mtn$ and $\alpha$ are smooth. Given $t\in T^n$ and $x\in \mtn$,
we abbreviate $\alpha(t,x)$ to $t\cdot x$ wherever possible. We assume
that $\alpha$ is locally equivalent to the standard action $\bC^n$, by
insisting that every point $x$ of \mtn\ lies in some neighbourhood $V$,
closed under the action of $\alpha$, for which there is a {\it
$\theta_x$-equivariant\/} diffeomorphism $h\colon V\rightarrow\bC^{n}$;
in other words,
\begin{equation}\label{locsta}
h(t\cdot y)=\theta_x(t)\cdot h(y)
\end{equation}
for some automorphism $\theta_x$ of $T^n$, and all $t\in T^n$ and
$y\in V$. Given a simple $n$-polyhedron $P^n$, we describe
\mtn\ as a {\it toric manifold over} $P^n$ whenever there exists a
smooth projection $\pi\colon\mtn\rightarrow P^n$ whose fibres are the
orbits of $\alpha$. We may display this information as the quadruple
$(\mtn,\alpha,\pi,P^n)$, and refer to $P^n$ as the {\it base
polyhedron}.

It is customary to insist that $\mtn$ and $P^n$ should be compact,
but the extra generality will prove helpful in Sections \ref{scs}
and \ref{cs} below.

Each face of $P^n$ of codimension $k$ is the image under $\pi$ of the
fixed point set of some $k$-dimensional subtorus, for all $0\leq k\leq
n$; for example, the vertices are the image of the fixed points, and the
boundary $\partial P^n$ is the image of the points on which $T^n$ fails
to act freely. The maps $h$ descend to local diffeomorphisms between
$P^n$ and the cone $\brpn$, yielding charts for $P^n$ as a manifold with
corners. In particular, the charts based on open subsets $U_v\subset
P^n$ correspond to a finite $T^n$-invariant atlas for $\mtn$, each of
whose open sets $V_x$ contains a single fixed point $x=\pi^{-1}(v)$. It
follows from \cite{daja:cpc} that $\pi$ admits smooth right inverses
$P^n\rightarrow\mtn$, from which we select a {\it preferred section}
$s$, transverse to the orbits. Any other choice differs from $s$ by some
map $P^n\rightarrow T^n$, and is therefore homotopic to $s$ through
right inverses because $P^n$ is contractible. We note that $P^n$ and the
quotient $\mtn/T^n$ are diffeomorphic as manifolds with corners.

Every facet $F$ of $P^n$ determines a subspace $\pi^{-1}(F)$,
readily seen to be a submanifold $\xft\subset\mtn$ with isotropy
subgroup a circle $T(F)$ in $T^n$. As $F$ ranges over
$\mathcal{F}_G$ for some face $G$ of codimension $k$, the
$X(F)^{2(n-1)}$ intersect transversally in a submanifold
$X(G)^{2(n-k)}$, whose isotropy subgroup $T(G)$ is a $k$-dimensional
subtorus, and is generated by the circles $T(F)$. We therefore have
a {\it characteristic map} $\lambda\colon
\mathfrak{L}_F(P)\rightarrow \mathfrak{L}_S(T^n)$ into the lattice
of subtori of $T^n$, which preserves the corresponding concept of
rank. In this way, we associate the {\it characteristic pair}
$(P^n,\lambda)$ to $(\mtn,\alpha,\pi,P^n)$.

Now let us reverse this process by starting with a pair
$(P^n,\lambda)$, where $\lambda$ is a rank-preserving map of the
lattices above. Note that each point $q$ of $\partial P^n$ lies in
the relative interior of a unique face $G(q)$. We use this data to
construct the identification space
\begin{equation}\label{dtm}
(T^n\times P^n)/\!\sim\,,
\end{equation}
where $(t,q)\sim (u,q)$ if and only if $tu^{-1}$ lies in the
subtorus $\lambda(G(q))$ of $T^n$. Multiplication on the first
coordinate defines an action of $T^n$ on the resulting space, with
quotient $P^n$. Whenever $q$ lies in the interior of $P^n$ the
equivalence classes $(t,q)$ are singletons, and have trivial
isotropy subgroups; at the other extreme, the fixed points consist
of the equivalence classes $(1,v)$, where $v$ ranges over the
vertices of $P^n$. Just as $P^n$ is covered by the open sets $U_v$,
based on the vertices and diffeomorphic to $\brpn$, so the
identification space is covered by open sets $(T^n\times
U_v)/\!\!\sim\,$, centred on the fixed points $(1,v)$ and
homeomorphic to $(T^n\times\brpn)/\!\!\approx$, and therefore to
$\bC^n$. With this structure, the rank-preserving properties of
$\lambda$ ensure that the identification space is a toric manifold
over $P^n$, which we say is {\it derived from\/} $(P^n,\lambda)$.

Given two toric manifolds over the same polyhedron $P^n$, we deem them
to be indestinguishable whenever they are linked by some
$\theta$-equivariant diffeomorphism (in the sense of \eqref{locsta})
which covers the identity map on $P^n$; here $\theta$ is an automorphism
of the torus $T^n$, and so induces an automorphism $\theta_*$ of the
lattice $\mathfrak{L}_S(T^n)$. Any such diffeomorphism descends to a
{\it $\theta$-translation\/} of characteristic pairs, in which the two
charactersitic maps differ by $\theta_*$. When $\theta$ is the identity,
these concepts reduce to equivariant diffeomorphism of toric manifolds
and equality of characteristic pairs, respectively. Two
$\theta$-equivariant diffeomorphisms $f$ and $f'$ are {\it equivalent}
whenever there exist equivariant diffeomorphisms $h_1$ and $h_2$ such
that $f\cdot h_1=h_2\cdot f'$.

\begin{prop}\label{oocop}
For any automorphism $\theta$, the assignment of characteristic pairs
defines a bijection between equivalence classes of $\theta$-equivariant
diffeomorphisms of toric manifolds, and $\theta$-translations of pairs
$(P^n,\lambda)$. 
\end{prop}
\begin{proof}
To prove bijectivity, we show that the inverse assignment is given by
taking derived toric manifolds; to each $\theta$-translation
$(P^n,\lambda)\rightarrow (P^n,\theta_*(\lambda)$ we associate the
$\theta$-equivariant diffeomorphism 
$\theta\times 1\colon(T^n\times P^n)/\!\!\sim\;\rightarrow (T^n\times
P^n)/\!\!\sim_\theta\,$, 
where $(t,q)\sim_\theta(u,q)$ if and only if
$tu^{-1}\in\theta_*(\lambda)(G(q))$. 

It follows directly from the definitions that $\theta\times 1$ descends
to the original $\theta$-translation $(P^n,\lambda)\rightarrow
(P^n,\theta_*(\lambda)$ of characteristic pairs. If, on the other hand,
we start with a $\theta$-equivariant diffeomorphism $f\colon
M_1^{2n}\rightarrow M_2^{2n}$ (or its equivalent), then $\theta\times 1$
is derived from the corresponding $\theta$-translation of characteristic
pairs. But the preferred section $s_1$ for $M_1^{2n}$ automatically
extends to an equivariant diffeomorphism $S_1\colon(T^n\times
P^n)/\!\!\sim\;\rightarrow M_1^{2n}$, and the section $s_2=f\cdot s_1$
extends to an eqivariant diffeomorphism $S_2\colon(T^n\times
P^n)/\!\!\sim_\theta\;\rightarrow M_2^{2n}$; thus $f\cdot
S_1=S_2\cdot(\theta\times 1)$, whence $f$ and $\theta\times 1$ are
equivalent, as required.
\end{proof}

In subsequent sections it will be convenient to replace
$(\mtn,\alpha,\pi,P^n)$ by its {\it derived form\/} \eqref{dtm}, and
use $S$ to transfer our constructions back to \mtn. We
abbreviate $(T^n\times P^n)/\!\!\sim\;$ to $\dfmtn$.

The {\it facial submanifolds} $X(G)^{2(n-k)}$ are central to
cobordism calculations, and form a lattice $\mathfrak{L}_X(\mtn)$
which is isomorphic to $\mathfrak{L}_F(P)$. We write $\nu(G)$ for
the normal $2k$-bundle of the embedding $X(G)^{2(n-k)}\subset\mtn$.
We may assume that $T(G)$ acts on the fibres of $\nu(G)$
isometrically with respect to a $T^n$-invariant metric; the
transformations acting tangentially form an $(n-k)$-dimensional
subtorus $T^\top(G)$, which splits $T^n$ as $T(G)\times
T^\top(G)$ and invests $X(G)^{2(n-k)}$ with its own toric
structure. We refer to this action as the {\it restriction\/} of
$\alpha$, and note that different choices of basis for $T^\top(G)$
correspond to $\theta$-equivariantly diffeomorphic versions of
$X(G)^{2(n-k)}$. Thus $\mathfrak{L}_X(\mtn)$ is a lattice of
subtoric manifolds.

Under $S$ the facial submanifold $\xgk$ corresponds to the
identification subspace $(T^\top(G)\times G^{n-k})/\!\!\sim\,$,
which we may equate with the derived form $\dfxgk$ once a basis is
chosen for $T^\top(G)$.

We are particularly interested in three families of toric manifolds.
They all happen to be toric varieties, but admit alternative stably
complex structures which arise naturally in the context of complex
cobordism theory. In giving their description, we denote a generic
element $t$ of $T^n$ by $(t_1,\dots,t_n)$, and for each $1\leq r\leq
n$ we write $T_r$ for the $r$th coordinate circle, defined by
$t_k=1$ unless $k=r$; we write the diagonal circle as $T_\delta$. In
each case we leave to readers the task of verifying that the action is
locally standard .

\begin{exa}\label{cpnex}
Complex projective space $\CPn$ is a toric manifold with respect to
the action induced by $t\cdot z=(t_1z_1,\dots,t_nz_n,z_{n+1})$ on the
unit sphere in $\bC^{n+1}$, and the projection
$\pi([z])=(|z_1|^2,\dots,|z_n|^2)$ onto the $n$-simplex $\splx{n}$.

In the notation of \eqref{nsplx}, the facial submanifold $X(D_r)$ is
a copy of $\CP{n-1}$; the normal bundle $\nu(D_r)$ is isomorphic (as
real $2$-plane bundles) to the Hopf bundle $\zeta(n-1)$, and is the
restriction of $\zeta(n)$. The characteristic map is given by
$\lambda(T^{D_r})=T_r$ for $1\leq r\leq n$, and
$\lambda(T^{D_{n+1}})=T_\delta$. Each $X(D_S)$ is a copy of
$\CP{n-|S|}$, and the lattice $\mathfrak{L}_X(\CPn)$ is Boolean of
rank $n$.

The stable tangent bundle arising from the complex algebraic
structure is given by a canonical isomorphism
$\tau(\CPn)\oplus\bC\cong(n+1)\oz(n)$ \cite{st:nct}.
\end{exa}
\begin{exa}\label{bnex}
The bounded flag manifold $B_n$ is described in \cite{bura:fml}, and
consists of complete flags $U$ in $\bC^{n+1}$ for which $U_k$
contains the subspace $\bC^{k-1}$ (spanned by the first $k-1$
standard basis vectors) for $2\leq k\leq n$; thus $U$ is equivalent
to a sequence of lines $L_k<\bC_k\oplus L_{k+1}$ for $1\leq k\leq
n$, where $\bC_k$ denotes the $k$th coordinate line, and
$L_{n+1}=\bC_{n+1}$. Then $B_n$ is a toric manifold with respect to
the action induced by $t\cdot z=(t_1z_1,\dots,t_nz_n,z_{n+1})$ on
$\bC^{n+1}$, and projection $\pi(U)=(\pi(L_1),\dots,\pi(L_n))$ onto
the $n$-cube $I^n$, where $\pi(L_k)$ is defined in
$\bR_{\scriptscriptstyle\geqslant}$ by projecting a unit vector onto
$\bC_k$ and taking the square of its modulus. For each $1\leq k\leq
n$, complex line bundles $\gamma_k(n)$ and $\rho_k(n)$ are defined
over $B_n$ by assigning to any bounded flag $U$ the line $L_k$ and
the orthogonal complement $L_{k,k+1}$ of $L_k$ in $\bC_k\oplus
L_{k+1}$ respectively. By convention we set $\gamma_0(n)$ and
$\gamma_{n+1}(n)$ to be trivial, and identify $\rho_0(n)$ with
$\gamma_1(n)$.

In the notation of \eqref{ncube}, the facial submanifold $X(C_r^0)$
is a copy of $B_{n-1}$, whose flags lie in
$\bC^{\{1,\dots,n+1\}\setminus r}$; the normal bundle $\nu(C_r^0)$ is
the restriction of $\gamma_r(n)$. On the other hand, $X(C_r^1)$ is a
copy of $B_{r-1}\times B_{n-r}$, where the flags of the factors lie
in $\bC^r$ and $\bC^{\{r+1,\dots,n+1\}}$ respectively; $\nu(C_r^1)$
is the restriction of $\rho_r(n)$. The submanifolds $X(C_r^0)$ and
$X(C_r^1)$ are labelled in \cite{bura:fml} as
$Y_{\{1,\dots,n\}\setminus r}$ and $X_{\{1,\dots,n\}\setminus r}$
respectively. The characteristic map is given by
$\lambda(T^{C_r^0})=T_r$ and $\lambda(T^{C_r^1})=T_\delta<T^r$,
where the latter is embedded in $T^n$ via the first $r$ coordinates,
for $1\leq r\leq n$. The lattice $\mathfrak{L}_X(B_n)$ is isomorphic
to $\mathfrak{L}_F(I^n)$.

Each $B_n$ is the sphere bundle of $\gamma_1\oplus\bR$ over
$B_{n-1}$. As detailed in \cite{ra:ocb}, this leads to a stably
complex structure
$\tau(B_n)\oplus\bR^2\cong\oplus_{k=2}^{n+1}\gamma_k(n)\oplus\bC$,
which plays an important r\^ole in complex cobordism theory despite
bounding the associated disk bundle.
\end{exa}
\begin{exa}\label{bijex}
The manifold $B_{i,j}$ (for integers $0\leq i\leq j$) is introduced
in \cite{bura:tmc}, and consists of pairs $(U,W)$, where $U$ is a
bounded flag in $\bC^{i+1}$ and $W$ is a line in
$U_1^\perp\oplus\bC^{j-i}$. So $B_{i,j}$ is a smooth
$\CP{j-1}$-bundle over $B_i$. It has dimension $2(i+j-1)$, and is a
toric manifold with respect to the action induced by
\[
t\cdot(z,w)=
(t_1z_1,\dots,t_iz_i,z_{i+1},t_{i+1}w_1,\dots,t_{i+j-1}w_{j-1},w_j)
\]
on $\bC^{i+1}\times(U_1^\perp\oplus\bC^{j-i})$, where the
coordinates of $w$ are chosen with respect to the decomposition
$L_{1,2}\oplus L_{2,3}\oplus\dots\oplus
L_{i,i+1}\oplus\bC^{\{i+2,\dots,j+1\}}$. Projection onto the product
$I^i\times\splx{j-1}$ is defined as $\pi(U,W)=(\pi(U),\pi(W))$, by
combining Examples \ref{cpnex} and \ref{bnex}. For each $1\leq k\leq
i$, complex line bundles $\gamma_k(i)$ and $\rho_k(i)$ are defined
over $B_{i,j}$ by pullback from $B_i$; similarly, $\zeta$ is defined
by considering $W$ as a line in $\bC^{j+1}$.

The facial submanifolds $X(E_r^0)$ and $X(E_s)$ are copies of
$B_{i-1,j}$ and $B_{i,j-1}$ respectively, for all $1\leq r\leq i$
and $i+1\leq s\leq j$; the corresponding normal bundles $\nu(E_r^0)$
and $\nu(E_s)$ are the restrictions of $\gamma_r(i)$ and $\zeta$.
The manifolds $X(E_r^1)$ and $X(E_s)$ for $1\leq s\leq i$ are new.
The characteristic map is given by
$\lambda(T^{C_r^0})=\{(t,t^{-1})\}$ in $T_r\times T_{i+r}$, and
$\lambda(T^{C_r^1})=\{(t,\dots,t,t^{-1},\dots,t^{-1})\}$ in
$T^r\times T^{r-1}$, where $T^r$ and $T^{r-1}$ are embedded in
$T^{i+j-1}$ by the first $r$ and the $\{i+1,\dots,i+r-1\}$
coordinates respectively; and also by $\lambda(T^{D_s})=T_{i+s}$ for
$1\leq s\leq j-1$ and $\lambda(T^{D_j})=T_\delta<T^{j-1}$, where
$T^{j-1}$ is embedded in $T^{i+j-1}$ as the last $j-1$ coordinates.
The lattice $\mathfrak{L}_X(B_{i,j})$ is isomorphic to the product
$\mathfrak{L}_F(I^i)\times\mathfrak{L}_F(\splx{j-1})$.

The projection onto $B_i$ and the classifying map of $\zeta$
together provide a smooth embedding of $B_{i,j}$ in
$B_i\times\CP{j}$, whose normal bundle is $\gamma_1(i)\otimes\zeta$.
Combining this with the bounding structure on $B_i$ and the varietal
structure on $\CP{j}$ yields the isomorphism
$\tau(B_{i,j})\oplus(\gamma_1(i)\otimes\oz)\oplus\bC\cong
\oplus_{k=2}^{i+1}\gamma_k(i)\oplus(j+1)\oz$ which defines the
stably complex structure used implicitly in \cite{bura:tmc}.
\end{exa}

Further examples are provided by taking connected sums of the above, as
outlined in \cite{daja:cpc}; however, the resulting tangent bundles are
rarely complex. As we now explain, the best we can generally expect is a
complex structure on the {\it stable\/} tangent bundle.

\section{Stably complex structure}\label{scs}

In order to describe our stably complex structures with appropriate
precision, we need to assign a collection of orientations to each
toric manifold. Davis and Januszkiewicz \cite{daja:cpc} sometimes
incorporate equivalent information into their notion of characteristic
map, but they do so implicitly, and without considering the dependence
of the resulting structures on their choice.

Given any toric manifold $(\mtn,\alpha,\pi,P^n)$ and any facet $F$ of
$P^n$, the action of $T(F)$ allows us to interpret the 2-plane bundle
$\nu(F)$ as a complex line bundle. Two complex structures are
possible, which differ by conjugation and correspond to opposite
orientations.  An {\it omniorientation\/} of $(\mtn,\alpha,\pi,P^n)$
consists of a choice of such orientation for every facet $F$; there
are therefore $2^m$ omniorientations in all, each of which is
preserved by $\alpha$ because $T^n$ is connected. By transversality,
an omniorientation determines an orientation (and also a complex
structure) for $\nu(G)$, given any face $G$ of $P^n$.

By analogy, we refer to a characteristic pair $(P^n,\lambda)$ as
{\it directed} if the circle $\lambda(F)$ is oriented for every
facet $F$ of $P^n$. We may then replace the lattice map $\lambda$ by
an epimorphism $\ell\colon T^\mathcal{F}\rightarrow T^n$, which
encodes each of the isomorphisms $T\rightarrow T(F)$ determined by
the orientation of the latter. We label $\ell$ a {\it directed\/}
characteristic map, or {\it dicharacteristic}, and write
$(P^n,\ell)$ for a directed characteristic pair; each $\lambda$ is
represented by $2^m$ distinct dicharacteristics.

The complex structures implicit in \cite{daja:cpc} need careful
interpretation precisely because $\lambda$ is used there to denote
both characteristic and dicharacteristic maps.

The characteristic pair of an omnioriented toric manifold is obviously
directed, and the toric manifold derived from a directed
characteristic pair is omnioriented. For any automorphism $\theta$ of
$T^n$, we insist that a $\theta$-equivariant diffeomorphism between
omnioriented toric manifolds should respect each of the $2^m$ facial
orientations; correspondingly, a $\theta$-translation of directed
pairs must satisfy $\ell_2=\theta\cdot\ell_1$. In this context, the
following extension of Proposition \ref{oocop} is immediate.
\begin{prop}\label{oocop2}
The assignment of directed characteristic pairs defines a bijection
between equivalence classes of $\theta$-equivariant diffeomorphisms of
omnioriented toric manifolds, and $\theta$-translations of pairs
$(P^n,\ell)$. 
\end{prop}
Transversality ensures that an omniorientation of $\mtn$ restricts to
an omniorientation of any facial submanifold $X(G)^{2(n-k)}$. If the
former corresponds to the dicharacteristic $\ell$ under Proposition
\ref{oocop2}, the latter corresponds to its restriction
\begin{equation}\label{dicharg}
\ell_{X(G)}\colon T^{\mathcal{F}(G)}\rightarrow T^\top(G)
\end{equation}
under the partition \eqref{partfacs} of $\mathcal{F}(P)$.

For each omniorientation of $(\mtn,\alpha,\pi,P^n)$ we now construct the
induced complex structures on $\mtn$. We focus initially on the base
polyhedron, which we assume to be defined in $\bR^n$ for convenience.

We recall the presentation of $P^n$ as a matrix inequality, and
interpret the $\mathcal{H}\times n$ matrix $A_P$ as a linear
transformation $A_P\colon\bR^n\rightarrow\bR^\mathcal{F}$. We abbreviate
the $n$-dimensional image $A_P(\bR^n)$ to $V_P$, and write $V_P^\perp$
for its $(m-n)$-dimensional orthogonal complement in $\bR^\mathcal{F}$
(with respect to the standard inner product). Since the points of $P^n$
are specified by the constraint $A_Px\geq b$, it follows that the
intersection of the affine subspace $V_P-b$ with the positive cone
$\brpf$ is a copy of $P^n$; it is embedded in $\bR^\mathcal{F}$ as the
space of functions $\{d(p,\;):\mathcal{F}\rightarrow\brp\}$, where
$d(p,F)$ is the euclidean distance between $p$ and the hyperplane
defining $F$ for each $p\in P^n$ and each facet $F$. We refer to this
embedding as $d_\mathcal{F}$. We sometimes identify $P^n$ with its image
$d_\mathcal{F}(P^n)$, which actually lies in the subspace
\begin{equation}\label{wrpdef}
\wrp=\{f\colon\mathcal{F}\rightarrow\brp
\;\text{such that}\;
f^{-1}(0)\in\mathfrak{L}_F(P)\}
\end{equation}
of $\brpf$.

For any subset $\mathcal{G}\subseteq\mathcal{F}$ of facets, we may
realise $\bR^\mathcal{G}$ as a subspace of $\bR^\mathcal{F}$ by
choosing $F$-coordinates to be $0$ for all $F$ in
$\mathcal{F}\setminus\mathcal{G}$. Thus $\wrp$ consists of the union
of open cones $\bigcup\bR_{\scriptscriptstyle >}^{\mathcal{C}_G}$,
where $\mathcal{C}_G$ denotes the complement of $\mathcal{F}_G$ in
$\mathcal{F}$, and $G$ ranges over $\mathfrak{L}_F(P)$; the union is
topologised by embedding the cones in $\brpf$ in the obvious
fashion, so that the interior of $G$ is embedded in
$\bR_{\scriptscriptstyle >}^{\mathcal{C}_G}$ for each face $G$.
Clearly $\wrp$ is a noncompact $m$-dimensional manifold with corners,
and never contains the zero vector.

The open cone $\brsf$ is an abelian topological group under
coordinatewise multiplication $*$. It decomposes as 
$\exp(V_P^\perp)\times\exp(V_P)$ by exponentiating the additive
splitting of $\bR^\mathcal{F}$ as $V_P^\perp\oplus V_P$. The group
$\exp(V_P^\perp)$ therefore acts smoothly on $\brsf$ by $*$, with
quotient space $\exp(V_P)$. This action restricts to each embedded cone
$\brsg$, and extends to $\brpg$ and $\wrp$.
\begin{prop}\label{strwrpn}
As manifold with corners, $\wrp$ is canonically diffeomorphic to the
cartesian product $P^n\times\exp(V_P^\perp)$.
\end{prop}
\begin{proof}
Given $p\in P^n$ and $a\in V_P^\perp$, we consider tangents to the
orbit $\exp(V_P^\perp)*p$ at $\exp(a)*p$. One such has direction
vector $a*\exp(a)*p$, whose inner product with $a$ is given by
$\sum_\mathcal{F}a_F^2\exp(a_F)p_F$. Since this quantity is strictly
positive the orbit meets $P^n$ only at $p$, (and the intersection is
transverse). On the other hand, given any point $x$ in the open cone
$\bR_{\scriptscriptstyle >}^{\mathcal{C}_G}$, the orbit
$\exp(V_P^\perp)*x$ meets $P^n$ in an interior point of $G$, for
each face $G$; this follows by taking logarithms and considering the
decomposition of $\bR^\mathcal{G}$ into $V_G^\perp\oplus V_G$. The
required diffeomorphism is therefore given by the map
$(p,\exp(a))\mapsto\exp(a)*p$.
\end{proof}
\begin{cor}\label{frembed}
The embedding of $P^n$ in $\brpf$ as manifolds with corners has trivial
normal bundle; each choice of basis for $V_P^\perp$ provides a framing. 
\end{cor}
\begin{proof}
Since exponentiation is a diffeomorphism, $\wrp$ is a tubular
neighbourhood of the embedding; each choice of basis for $V_P^\perp$
therefore trivialises the normal bundle.
\end{proof}

We consider the identification space
$(T^\mathcal{F}\times\wrp)/\!\!\approx$, denoted by $W(P)$, as a
complexified form of the tubular neighbourhood. Such a space has also
been introduced by Buchstaber and Panov \cite{bupa:tac}, as an
extension of a construction for toric varieties \cite{ba:qcr},
\cite{co:rdt}. By \eqref{wrpdef}, $W(P)$ embeds in $\bC^\mathcal{F}$
as the space of complex valued functions whose zero-set is
$\mathcal{F}_G$ for some face $G$ of $P^n$. The multiplicative group
$(\bC_\times)^\mathcal{F}$ of vectors with nonzero coordinates acts on
$\bC^\mathcal{F}$ by $*$, and the subgroups $T^\mathcal{F}$ and
$\exp(V_P^\perp)$ restrict to $W(P)$ by construction.

We now turn to the omniorientation of $(\mtn,\alpha,\pi,P^n)$, and
write $K(\ell)$ for the kernel of the dicharacteristic; it is an
$(m-n)$-dimensional subtorus of $T^\mathcal{F}$, and therefore also
acts on $W(P)$ by $*$. The quotient of $T^\mathcal{F}$ by $K(\ell)$
is, by definition, the original torus $T^n$, from which we deduce
that the projection 
\begin{equation}\label{kexpbun}
(T^\mathcal{F}\times\wrp)/\approx\;\;\longrightarrow\; (T^n\times
P^n)/\sim
\end{equation}
displays $W(P)$ as a smooth, principal
$K(\ell)\times\exp(V_P^\perp)$-bundle over the derived form
$\dfmtn$. We abbreviate $K(\ell)\times\exp(V_P^\perp)$ to $H(\ell)$;
since it is a subgroup of $(\bC_\times)^\mathcal{F}$, the embedding
of $W(P)$ in $\bC^\mathcal{F}$ is $H(\ell)$-equivariant.

The tangent bundle of $W(P)$ inherits a natural complex structure
from that of $(\bC_\times)^\mathcal{F}$, and its quotient by the
action of $H(\ell)$ provides our stably complex structure on $\mtn$. The
details, however, need care; they involve extending Sczarba's analysis
of \cite{sz:tbf} to \eqref{kexpbun}, circumventing his restriction to
compact fibres. We obtain a canonical isomorphism 
\begin{equation}\label{scziso}
\tau(\dfmtn)\oplus\tau^\|\cong\sigma(\mathcal{F})
\end{equation}
of real $2m$-plane bundles, where $\tau^\|$ is the quotient of the
tangents along the fibres by $H(\ell)$, and $\sigma(\mathcal{F})$
is the $\bC^\mathcal{F}$-bundle associated to \eqref{kexpbun}. We
equip each of these bundles with the standard inner product, and
insist that $H(\ell)$ acts on the fibres $\bC^\mathcal{F}$
by projection onto its maximal compact subgroup $K(\ell)$. Of
course $\sigma(\mathcal{F})$ is isomorphic to the $m$-fold sum of
complex line bundles $\oplus_\mathcal{F}\sigma(F)$, where $\sigma(F)$
has total space $W(P)\times_{H(\ell)}\bC^F$.

\begin{thm}\label{cplxstr}
An omniorientation of a toric manifold induces a stably complex
structure on the derived form; it is defined uniquely up to homotopy.
\end{thm}
\begin{proof}
Following \eqref{scziso}, we must identify $\tau^\|$ with
$\bR^{2(m-n)}$. Any choice of basis for the Lie algebra of $H(\ell)$
will have this effect; since the space of bases has two connected
components, it actually suffices to give an orientation. But $H(\ell)$
is the kernel of an epimorphism $\bC^\mathcal{F}_\times\rightarrow 
T^n\times\exp(V_P)$ and $\bC^\mathcal{F}$ is canonically oriented, so
it remains to orient the domain; since the latter is isomorphic to
$(\bC_\times)^n$, we may simply apply the standard orientation of
$T^n$. The resulting isomorphism 
\begin{equation}\label{scsdf} 
\tau(\dfmtn)\oplus\bR^{2(m-n)}\cong\sigma(\mathcal{F})
\end{equation}
provides the structure we seek.
\end{proof}
If we use the opposite orientation for $T^n$, and hence for $T^n\times
\exp(V_P)$, we obtain a second stably complex structure. This is
compatible with the opposite orientation on $\dfmtn$, and its complex
cobordism class is the negative of that represented by \eqref{scsdf}. We
emphasise that the isomorphism \eqref{scsdf} does not depend on any
ordering of $\mathcal{F}$.

Theorem \ref{cplxstr} gives a global description for any toric
manifold, with its induced complex structure, as the quotient of a
complex space; given a nonsingular toric variety \cite{co:rdt}, it
yields the stabilisation of the underlying complex structure.

To continue our investigation of induced complex structures, we need a
technical lemma. It considers directed characteristic pairs
$(P_1^n,\ell_1)$ and $(P_2^n,\ell_2)$, with omnioriented derived forms
$M_1^{2n}$ and $M_2^{2n}$ respectively, and assumes given closed
halfspaces $H_1$ and $H_2$ in $\bR^n$, and a diffeomorphism $f\colon
P_1^n\setminus H_1\rightarrow P_2^n\setminus H_2$ as manifolds with
corners. We abbreviate $P_1^n\setminus H_1$ to $O^n$, and write
$\mathcal{C}_1$ and $\mathcal{C}_2$ for the sets of facets contained
in $H_1$ and $H_2$ respectively. We partition $\mathcal{F}_1$ as
$\mathcal{E}\cup\mathcal{C}_1$ and $\mathcal{F}_2$ as
$\mathcal{E}\cup\mathcal{C}_2$, where $\mathcal{E}$ consists of those
facets which intersect $O^n$; we use $f$ to identify $\mathcal{E}$
with the set of facets intersecting $f(O^n)$.

\begin{lem}\label{techlem}
If $f$ preserves dicharacteristics on $\mathcal{E}$, it lifts to a
an equivariant diffeomorphism
\[
f^+\colon(T^n\times O^n)/\!\sim\;\;\longrightarrow(T^n\times
f(O^n))/\!\sim\;\;;
\]
$f^+$ respects the stably complex structures obtained by
restriction from those induced on $M_1^{2n}$ and $M_2^{2n}$ by their
respective omniorientations.
\end{lem}
\begin{proof}
The existence of $f^+$ is assured by the fact that $\ell_1$
and $\ell_2$ agree on $\mathcal{E}$; we write their common restriction
as $\ell$.

Our data implies that $d_{\mathcal{F}_1}(O^n)$ maps
diffeomorphically onto $d_\mathcal{E}(O^n)$ under the projection of
$\bR^{\mathcal{F}_1}$ onto $\bR^{\mathcal{E}}$, so that
$\exp(V_{P_1}^\perp)$ is isomorphic to
$\exp(V_O^\perp)\times\bR_{\scriptscriptstyle >}^{\mathcal{C}_1}$ in
$\bR_{\scriptscriptstyle\geqslant}^\mathcal{E}\times
\bR_{\scriptscriptstyle\geqslant}^{\mathcal{C}_1}$. Since 
$W^\bR(P_1)|_O$ is given by $\exp(V_{P_1}^\perp)*O^n$ in
$\bR_{\scriptscriptstyle\geqslant}^{\mathcal{F}_1}$, it is
equivariantly diffeomorphic to
$\exp(V_O^\perp)*O^n\times\bR_{\scriptscriptstyle
>}^{\mathcal{C}_1}$ with respect to the splitting of
$\exp(V_{P_1}^\perp)$; we write $W^\bR(O)$ for the subspace
$\exp(V_O^\perp)*O^n\subset
\bR_{\scriptscriptstyle\geqslant}^\mathcal{E}$. We also note that
$K(\ell_1)$ splits as $K(\ell)\times T^{\mathcal{C}_1}$. The stably
complex structure on $(T^n\times O^n)/\!\sim\;$ given by factoring
out the action of $H(\ell_1)$ on $W(P)|_O$ therefore differs from
that given by factoring out the action of $H(\ell)$ on $W(O)$ only
by the trivial summand $\bC^{\,\mathcal{C}_1}$, so that we may  
consider the latter as obtained from $M_1^{2n}$ by restriction.
Similar remarks apply to $(T^n\times f(O^n))/\!\sim\;$ and
$M_2^{2n}$.

To show that $f^+$ respects the restricted complex structures, we
then choose an isomorphism
$\exp(V_O^\perp)\rightarrow\exp(V_{f(O)}^\perp)$; this immediately
extends to an equivariant diffeomorphism $W^\bR(O)\rightarrow
W^\bR(f(O))$, and therefore to an equivariant diffeomorphism
$W(O)\rightarrow W(f(O))$ which preserves the action of
$T^\mathcal{E}$. The differential of the second diffeomorphism is
complex linear, and reduces to $df^+$ on the quotient tangent
bundles, as required.
\end{proof}

Our first corollary to Lemma \ref{techlem} deals with the reliance of
Theorem \ref{cplxstr} on the hyperplanes defining $P^n$, in apparent
contradiction to the fact that an omniorientation of $\mtn$ involves
only the action $\alpha$.
\begin{cor}\label{uniquecs}
Given a common omniorientation for $(\mtn,\alpha,\pi_1,P_1^n)$ and
$(\mtn,\alpha,\pi_2,P_2^n)$, the derived forms are equivariantly
diffeomorphic; the induced stably complex structure therefore depends
only on the combinatorial type of the base polyhedron.
\end{cor}
\begin{proof}
The data yields a diffeomorphism $f\colon P_1^n\rightarrow P_2^n$,
with $f\cdot\pi_1=\pi_2$; we then apply Lemma \ref{techlem} to $f$ with
$O^n=P_1^n$, and to $f^{-1}$ with $O^n=P_2^n$.
\end{proof}

Our second application relates two stably complex structures which are
naturally prescribed on $\dfxgk$ by Theorem \ref{cplxstr}, for any
codimension $k$ face $G$. One is induced by the restricted
omniorientation associated with the dicharacteristic $\ell_{X(G)}$ of
\eqref{dicharg}; the other is the restriction to $\dfxgk$ of the
structure induced on $\dfmtn$, using the complex structure given on
$\nu(G)$ by the omniorientation. 

We confirm that these are equivalent in Theorem \ref{restrcs}. The proof
involves an auxiliary polyhedron $R^n$, which is defined by expressing
$G^{n-k}$ as an intersection of halfspaces in $\bR^{n-k}$ and taking
products with $\bR_{\scriptscriptstyle\geqslant}^{\mathcal{F}_G}$; the
result is a simple $n$-polyhedron in
$\bR^{n-k}\times\bR^{\mathcal{F}_G}$, whose facets $\mathcal{F}(R)$ may
be partitioned as $\mathcal{F}(G)\cup\mathcal{F}_G$. Then we have
\begin{equation}\label{embr}
d_{\mathcal{F}(R)}(R^n)=d_{\mathcal{F}(G)}(G^{n-k})\times
\bR_{\scriptscriptstyle\geqslant}^{\mathcal{F}_G}
\quad\text{in}\quad\bR_{\scriptscriptstyle\geqslant}^{\mathcal{F}(R)}. 
\end{equation}
The restriction of $\ell$ to $\mathcal{F}(R)$ agrees with
$\ell_{X(G)}$ on $\mathcal{F}(G)$ and defines an omnioriented toric
manifold $L^{2n}$ over $R^n$. We invest $L^{2n}$ with the induced
stably complex structure, and note that $\dfxgk$ is a facial
submanifold, to which the omniorientations of $\dfmtn$ and $L^{2n}$ have
common restriction. 

\begin{thm}\label{restrcs}
The two stably complex structures on $\dfxgk$ are homotopic; that is, 
restriction and induction commute for $\dfmtn$.
\end{thm}
\begin{proof}
Considering $G$ as a face of $R^n$, we note that $W^\bR(R)=
W^\bR(G)\times\bR_{\scriptscriptstyle\geqslant}^{\mathcal{F}_G}$ in
$\bR_{\scriptscriptstyle\geqslant}^{\mathcal{F}(R)}$, that
$K(\ell_L)=K(\ell_{X(G)_\bullet})\times 1$ in
$T^{\mathcal{F}(G)}\times T^{\mathcal{F}_G}$, and that 
$V_{R}^\perp=V_G^\perp\times 0$ in
$\bR^{\mathcal{F}(G)}\times\bR^{\mathcal{F}_G}$ by \eqref{embr}.
Thus $W(R)=W(G)\times\bC^{\mathcal{F}_G}$, and $H(\ell_L)$ acts as
$H(\ell_{X(G)_\bullet})\times 1$. It follows that the normal bundle
$\nu$ of $\dfxgk$ in $L^{2n}$ has total space
$W(G)\times_{H(\ell_L)}\bC^{\mathcal{F}_G}$, and therefore that
restriction and induction commute for $L^{2n}$. 

Considering $G$ as a face of $P^n$, we may no longer appeal to
\eqref{embr}. Instead, we apply Lemma \ref{techlem} with $P_1^n=R^n$
and $P_2^n=P^n$; we let $H_2$ have complement an open tubular
neighbourhood $N(G)$ of $G$, and $f\colon R^n\rightarrow N(G)$ be any
diffeomorphism extending the identity on $G$. The Lemma provides an
equivariant diffeomorphism $f^+$ between $L^{2n}$ and an open tubular
neighbourhood of $\dfxgk$ in $\dfmtn$. By construction, $f^+$ is
compatible with the stable complex structures induced by $\ell_L$ and
$\ell$ respectively, and defines an isomorphism between $\nu$ and the
normal bundle $\nu_\bullet(G)$ of $\dfxgk$ in $\dfmtn$. This
isomorphism is well-defined up to homotopy, and therefore confirms
that restriction and induction commute for $\dfmtn$.
\end{proof}

Before pulling our constructions back to $\mtn$ we consider how the
bundles $\sigma(F)$ restrict to $\dfxgk$, for any facet $F$ and any
face $G$ of codimension $k$.  
\begin{prop}\label{propsigg}
For any facet $D$ disjoint from $G$, the restriction
$\sigma(D)|_{X(G)_\bullet}$ is trivial; on the other hand, 
$\sigma(\mathcal{F}_G)|_{X(G)_\bullet}$ is isomorphic to
$\nu_\bullet(G)$.
\end{prop}
\begin{proof}
The first statement follows from the proof of Lemma \ref{techlem} 
by choosing $P_1^n=P^n$, and letting $O^n$ be a tubular neighbourhood
of $G$; then $D$ lies in $\mathcal{C}_1$, and
$\sigma(D)|_{X(G)_\bullet}$ is the corresponding coordinate line
bundle in the trivial summand $\bC^{\,\mathcal{C}_1}$. For the second
statement, we note that the proof of Theorem \ref{restrcs} identifies
the total space of $\nu_\bullet(G)$ with
$W(P)|_{X(G)_\bullet}\times_{H(\ell)}\bC^{\mathcal{F}_G}$. 
\end{proof}

Proposition \ref{propsigg} leads to an alternative description of
$(S^{-1})^*\sigma(F)$, which simplifies subsequent calculations in
cobordism theory. We express the orientation of $\nu(F)$ as an
integral Thom class in the cohomology group $H^2(M(\nu(F)))$,
represented by a complex line bundle over the Thom complex
$M(\nu(F))$. We pull this back along the Pontryagin-Thom collapse
$\mtn\rightarrow M(\nu(F))$, and label the resulting bundle $\rho(F)$.
\begin{lem}\label{roisnorm}
The line bundles $\sigma(F)$ and $S^*\rho(F)$ are isomorphic over
$\dfmtn$.  
\end{lem}
\begin{proof}
Since $S$ restricts to a preferred section for $\xft$, so $S^{-1}$
pulls $\nu_\bullet(F)$ back to $\nu(F)$. From Proposition
\ref{propsigg}, we deduce that $(S^{-1})^*\sigma(F)$ is isomorphic to
$\nu(F)$ over $\xft$, and is trivial over the complement. But these
properties characterise $\rho(F)$.
\end{proof}

We refer to the $\rho(F)$ as the {\it facial bundles\/} of
$\mtn$, to distinguish them from the canonical line
bundles $L_F$ of algebraic geometry, defined when $\mtn$ is
also a toric variety. In fact $\rho(F)$ and $L_F$ are either
isomorphic or complex conjugate.

\begin{thm}\label{bigthm}
An omniorientation of a toric manifold $(\mtn,\alpha,\pi,P^n)$
induces a canonical stably complex structure on $\mtn$, which is
preserved by the action $\alpha$; for each facial submanifold
$\xgk$, the restriction of this structure is homotopic to that
induced by the restricted omniorientation.
\end{thm}
\begin{proof}
Pulling \eqref{scsdf} back along $S^{-1}$ yields the complex
structure on $\mtn$, which Lemma \ref{roisnorm} converts to an
isomorphism $\tau(\mtn)\oplus\bR^{2(m-n)}\cong\rho(\mathcal{F})$.
Different choices of preferred section $s$ yield homotopic
isomorphisms, and therefore homotopic stably complex structures,
because the corresponding diffeomorphisms $S$ are isotopic. By
Proposition \ref{propsigg}, the complex structure on $\nu(G)$ is
given by an isomorphism $\nu(G)\cong\rho(\mathcal{F}_G)$, which is
defined uniquely up to homotopy. The structures restrict to $\xgk$ as
claimed, by appeal to Theorem \ref{restrcs}.  
\end{proof}

Theorem \ref{bigthm} allows us to compute the complex bordism and
cobordism of an arbitrary toric manifold, as explained in Section
\ref{cc} below. It also shows that any choice of omniorientation leads
to an {\it equivariant\/} complex cobordism class, as defined, for
example, in Comeza\~na \cite{co:cce}. Equivariant complex cobordism is
currently under active development, and we expect the r\^ole played by
toric manifolds to be clarified in future work.

\section{Examples}\label{exs}

We consolidate our results by returning to the examples of Section
\ref{toma}. It is convenient to follow Davis and Januszkiewicz by
simplifying the dicharacteristic to a function which assigns to each
facet $F$ of $P^n$ a primitive vector $\ell(F)$ in $\bZ^n$; this is
obtained by applying the induced map $d\ell$ of Lie algebras to the
positively oriented unit tangent vector of $T^F$.

\begin{exa}\label{cpnext}
For $\CPn$ as in Example \ref{cpnex}, we note that $W^\bR(\splx{n})$
is $\bR_{\scriptscriptstyle\geqslant}^{n+1}\setminus 0$, and therefore
that $W(\splx{n})$ is isomorphic to $\bC^{n+1}\setminus 0$, by ordering
the facets of $\splx{n}$ as in \eqref{nsplx}. The dicharacteristic
chosen by Davis and Januszkiewicz, albeit implicitly, is
\begin{equation}\label{cpndich1}
\ell'(D_r)=
\begin{cases}
(0,\dots,0,1,0,\dots,0)&\text{for $1\leq r\leq n$}\\
(1,1,\dots,1)&\text{for $r=n+1$},
\end{cases}
\end{equation}
so that $K(\ell')$ is the subcircle $\{(t,\dots,t,t^{-1})\}$ of
$T^{n+1}$. Since the normal vector to $\splx{n}$ in
$\bR_{\scriptscriptstyle\geqslant}^{n+1}$ is 
$(1,\dots,1,n^{1/2})$, the action of $H(\ell')$ on $W(\splx{n})$ is
equivariantly diffeomorphic to that of $\bC_\times$ on
$\bC^{n+1}\setminus 0$ by
\[
x\cdot(z_1,\dots,z_n,z_{n+1})=(xz_1,\dots,xz_n,\overline{x}z_{n+1}).
\]
The composition of $S$ with the projection \eqref{kexpbun}
maps $(z_1,\dots,z_n,z_{n+1})$ to the line
$[z_1,\dots,z_n,\overline{z}_{n+1}]$, ensuring that the facial bundles
$\rho(D_r)$ are given by $\oz(n)$ for $1\leq r\leq n$, and $\zeta(n)$
for $r=n+1$. The omniorientation corresponding to \eqref{cpndich1}
induces a stably complex structure
$\tau(\CPn)\oplus\bR^2\cong n\,\oz(n)\oplus\zeta(n)$.

A second dicharacteristic $\ell$ arises by setting
$\ell(D_{n+1})=(-1,-1,\dots,-1)$; the corresponding omniorientation
induces $\tau(\CPn)\oplus\bC^1\cong(n+1)\oz(n)$, which stabilises the
structure of $\CPn$ considered as an algebraic variety. Both $\ell$
and $\ell'$ represent the characteristic map of Example \ref{cpnex}.
When $n=1$, the structure induced by $\ell$ represents a generator of
$\varOmega_2^U$, whereas that induced by $\ell'$ extends over the
$3$-disk, and represents zero.
\end{exa}

\begin{exa}\label{bnext}
For $B_n$ as in Example \ref{bnex}, we note that $W^\bR(I^n)$ is
$(\bR_{\scriptscriptstyle\geqslant}^{\{0,1\}}\setminus 0)^n$, and
therefore that $W(I^n)$ is isomorphic to $(\bC^2\setminus 0)^n$. The
characteristic map of Example \ref{bnex} is represented by the
dicharacteristic
\begin{equation}\label{bndich}
\ell(C^\epsilon_r)=
\begin{cases}
(0,\dots,0,-1,0,\dots,0)&\text{if $\epsilon=0$}\\
(-1,\dots,-1,0,\dots,0)&\text{if $\epsilon=1$},
\end{cases}
\end{equation}
for all $1\leq r\leq n$ (where the nonzero elements are in positions
$r$ and $1,\dots,r$ respectively). Thus $K(\ell)$ is the
$n$-dimensional subtorus
\[
\{(t_1,t_1^{-1}t_2,\dots,t_r,t_r^{-1}t_{r+1},\dots,
t_{n-1},t_{n-1}^{-1}t_n,t_n,t_n^{-1})\}
\]
of $T^{2n}$. Since the normal space to $I^n$ in
$\bR_{\scriptscriptstyle\geqslant}^{2n}$ is spanned by the $n$ vectors
$(0,\dots,0,1,1,0,\dots,0)$ (nonzero only in positions $2r$ and
$2r+1$), the action of $H(\ell)$ on $W(I^n)$ is equivariantly
diffeomorphic to that of $(\bC_\times)^n$ on $(\bC^2\setminus 0)^n$ by
\[
(x_1,\dots,x_n)\cdot(z_1,w_1,\dots,z_n,w_n)=
(x_1z_1,\overline{x}_1x_2w_1,\dots,x_nz_n,\overline{x}_nw_n).
\]
The composition of $S$ with the projection \eqref{kexpbun}
maps $(z_1,w_1,\dots,z_n,w_n)$ to the bounded flag for which
$L_r$ is spanned by the unit vector $l_r$, given by
$\overline{z}_r\overline{w}_r\dots\overline{w}_ne_r+\lambda l_{r+1}$
(where $e_r$ is the $r$th basis vector, $\lambda$ is the normalising
factor, and $l_{n+1}=e_{n+1}$), for each $1\leq r\leq n+1$. Since
$H(\ell)$ acts on $\bC^{C^0_r}$ by multiplication by $x_r$, and on
$\bC^{C^1_r}$ by multiplication by  $\overline{x}_rx_{r+1}$, the
associated facial bundles $\rho(C^0_r)$ and $\rho(C^1_r)$ are given by
$\gamma_r(n)$ and $\rho_r(n)$ respectively. The omniorientation
corresponding to \eqref{bndich} therefore induces a stably complex
structure
$\tau(B_n)\oplus\bR^{2n}\cong
\oplus_{r=1}^n(\gamma_r(n)\oplus\rho_r(n))$. When combined with
the canonical trivialisation of $\gamma_1(n)\oplus_{r=1}^n\rho_r(n)$,
this reduces to the bounding structure of Example \ref{bnex}.
\end{exa}

\begin{exa}\label{bijext}
For $B_{i,j}$ as in Example \ref{bijex}, we note that
$W(I^i\times\splx{j-1})$ is isomorphic to
$(\bC^2\setminus 0)^i\times\bC^j\setminus 0$. The characteristic
map of Example \ref{bijex} is represented by the dicharacteristic
\begin{equation}\label{bijdich}
\ell(E^\epsilon_r)=
\begin{cases}
(0,\dots,0,-1,0,\dots,0,1,0,\dots,0)&\text{if $\epsilon=0$}\\
(-1,\dots,-1,-1,0,\dots,0,1,\dots,1,0,\dots,0)&\text{if $\epsilon=1$}
\end{cases}
\end{equation}
for all $1\leq r\leq n$ (where the nonzero elements are in positions
$r$ and $i+r$, and $1,\dots,r$ and $i+1,\dots,i+r-1$ respectively),
and
\[
\ell(E_s)=
\begin{cases}
(0,\dots,0,1,0,\dots,0)&\text{for $1\leq s\leq j-1$}\\
(0,\dots,0,-1,\dots,-1)&\text{for $s=j$}
\end{cases}
\]
(where the nonzero elements are in positions $i+s$ and
$i+1,\dots,i+j-1$ respectively). Thus $K(\ell)$ is the
$(i+1)$-dimensional subtorus
\[
\{(t_1,t_1^{-1}t_2,\dots,t_{i-1},t_{i-1}^{-1}t_i,t_i,t_i^{-1},
t_1^{-1}t_2t,\dots,t_{i-1}^{-1}t_it,t_i^{-1}t,t,\dots,t)\}
\]
of $T^{2i+j}$. Combining Examples \ref{cpnext} and \ref{bnext},
we deduce that the associated facial bundles $\rho(E^0_r)$ and
$\rho(E^1_r)$ are given by $\gamma_r(i)$ and $\rho_r(i)\otimes\oz$
respectively, for each $1\leq r\leq i$, and that $\rho(E_s)$ is given
by $\oz$ for each $1\leq s\leq j$. The omniorientation
corresponding to \eqref{bijdich} therefore induces a stably complex
structure
\[
\tau(B_{i,j})\oplus\bR^{2(i+1)}\cong
\oplus_{r=1}^i\left(\gamma_r(i)\oplus(\rho_r(i)\otimes\oz)\right)
\oplus j\oz\,.
\]
Adding $\gamma_1(i)\otimes\oz$ to both sides and applying the
canonical trivialisation of $\gamma_1(i)\oplus_{r=1}^i\rho_r(i)$, we
obtain the stably complex structure of Example \ref{bijex} and
\cite{bura:tmc}.
\end{exa}

We also wish to consider products of these examples, and find it
equally convenient to discuss the general case. We assume given
omnioriented toric manifolds $(\mtn,\alpha,\pi,P^n)$ and
$(N^{2n},\beta,\mu,Q^n)$, with corrresponding dicharacteristic maps
$\ell_M$ and $\ell_N$. The facets of $P^n\times Q^n$ are of the form
$E\times Q^n$ and $P^n\times F$, where $E$ and $F$ range over
$\mathcal{F}(P)$ and $\mathcal{F}(Q)$ respectively, so we may define a
product dicharacteristic $\ell_{M\times N}\colon\mathcal{F}(P\times
Q)\rightarrow T^n\times T^n$ by assigning values $\ell_M(E)\times 1$ to
$E\times Q^n$ and $1\times\ell_N(F)$ to $P^n\times F$. This
corresponds to the product omniorientation of $(\mtn\times
N^{2n},\alpha\times\beta,\pi\times\mu,P^n\times Q^n)$.
\begin{prop}\label{prodcs}
The stably complex structure induced on $\mtn\times N^{2n}$ by the
product omniorientation is homotopic to the product of the structures
induced by the omniorientations of $\mtn$ and $N^{2n}$.
\end{prop}
\begin{proof}
By definition, there is a canonical diffeomorphism between
$W(P\times Q)$ and $W(P)\times W(Q)$, and a canonical isomorphism
between $H(\ell_{M\times N})$ and $H(\ell_M)\times H(\ell_N)$ which
preserves their respective actions on
$\bC^{\mathcal{F}(P)\sqcup\mathcal{F}(Q)}$ and
$\bC^{\mathcal{F}(P)}\times\bC^{\mathcal{F}(Q)}$. We obtain a
diffeomorphism 
$(\mtn\times N^{2n})_\bullet\rightarrow\dfmtn\times N^{2n}_\bullet$
which repects the induced stably complex structures, and the result
therefore follows by pullback along inverse preferred sections, as in
Theorem \ref{bigthm}. 
\end{proof}

In \cite{bura:tmc}, we showed that the $B_{i,j}$ are multiplicative
generators for the complex cobordism ring $\osu$ when
invested with the stably complex structure of Example
\ref{bijext}. So every $2n$-dimensional complex cobordism class may
be represented by a disjoint union of products
\begin{equation}\label{dup}
B_{i(1),j(1)}\times B_{i(2),j(2)}\times\dots\times B_{i(t),j(t)},
\end{equation}
where $\sum_{k=1}^t(i(k)+j(k))-2t=n$. Each such component is a toric
manifold with the product toric structure. This result is the
substance of \cite{bura:tmc} and may now be enriched by combining
Example \ref{bijext} with Proposition \ref{prodcs} to confirm that
the stably complex structures in question are induced by
omniorientations, and are therefore also preserved by the torus
action.

To give genuinely toric representatives (which are, by definition,
connected) for each cobordism class of dimension $>2$, it remains only
to replace the disjoint union of products \eqref{dup} with their
connected sum. This we do in Section \ref{cs} below.

\section{Cobordism calculations}\label{cc}

We now outline certain consequences of the results of Section
\ref{scs}. Our aim is to adapt Davis and Januszkiewicz's programme for
computing the integral homology and cohomology of a compact toric
manifold $(\mtn,\alpha,\pi,P^n)$, so as to apply directly to the complex
bordism groups $\osu(\mtn)$ and the complex cobordism ring
$\ous(\mtn)$. As always, we work with a fixed omniorientation.
We begin by summarising a few prerequisites concerning the bordism and
cobodism groups of a CW complex $X$ of finite type. These may be
found, for example, in the books of Stong \cite{st:nct} and Switzer
\cite{sw:ath}.

When the cells of $X$ lie only in even dimensions the
Atiyah-Hirzebruch spectral sequence collapses, and confirms that
generators for the free $\osu$-module $\osu(X)$ are given by the
bordism classes of any set of singular, stably complex manifolds whose
top cells correspond to the cells of $X$. The groups $\ous(X)$ are
then the $\Hom_{\varOmega^U_*}$-duals; when $X$ is itself a stably complex
manifold, the multiplicative structure of $\ous(X)$ may be extracted
from the intersection theory of the generating set by Poincar\'e
duality.

With these considerations in mind, we follow the opening gambit of
\cite{daja:cpc} by constructing a cell decomposition for $\mtn$. This
depends on choosing a generic direction in the ambient $\bR^n$, and so
determining an orientation for each edge of the $1$-skeleton of the
polytope $P^n$; the $1$-skeleton becomes a directed graph $Di(P)$. Any
vertex $v$ has indegree $m(v)$ and outdegree $n-m(v)$ for some integer
$0\leq m(v)\leq n$, and the $m(v)$ inward edges define an
$m(v)$-dimensional face $G_v$ of $P^n$. We write $\widehat{G}_v$ for the
subspace obtained by deleting all faces of $G_v$ disjoint from $v$,
which is therefore diffeomorphic to
$\bR_{\scriptscriptstyle\geqslant}^{m(v)}$, and write $e_v\subset\mtn$
for the subspace $\pi^{-1}(\widehat{G}_v)$. Since $e_v$ may be
identified with $\bC^{m(v)}$ it is a $2m(v)$-dimensional cell, and lies
within the open set $U_v$ of Section \ref{toma}; in fact $e_v$ and $U_v$
coincide precisely when $v$ is the sink of $Di(P^n)$, in which case
$e_v$ has dimension $2n$. The resulting decomposition of $\mtn$ has one
cell for each vertex of $P^n$, and all the cells are even
dimensional. The closure of $e_v$ is the facial submanifold
$X(G_v)^{2m(v)}$, which inherits the stably complex structure of Theorem
\ref{bigthm}.

\begin{prop}\label{osumtn}
The $\osu$-module $\osu(\mtn)$ is generated by the inclusions of the
facial submanifolds $\xgk$; none of these is null-cobordant, but they
are subject to nontrivial linear relations.
\end{prop}
\begin{proof}
The first statement follows from our introductory remarks, in view of
the cell decomposition defined above. Since there are $m$ submanifolds
$\xft$, but only $m-n$ two-cells in the decomposition, Poincar\'e
duality shows that there are $n$ linear relations amongst the
cobordism classes of the inclusions.   
\end{proof}
Proposition \ref{osumtn} highlights the remarkable fact that
$\osu(\mtn)$ is spanned by {\it embedded\/} submanifolds, each of
which is equipped with the restricted stably complex structure and is
itself a toric manifold.

The omniorientation determines $m$ cobordism Chern classes
$c_1(\rho(F))$, each lying in $\varOmega^2_U(\mtn)$ and Poincar\'e
dual to the inclusion $\xft\subset\mtn$ by construction of
$\rho(F)$. By transversality, any product $c_1(\rho(F_1))\cdots
c_1(\rho(F_k))$ is Poincar\'e dual to the inclusion of the facial
submanifold $X(F_1\cap\dots\cap F_k)$; if the intersection of the
facets is empty, the bordism and cobordism classes vanish
together. So the lattice $\mathfrak{L}_X(\mtn)$ maps into both
$\osu(\mtn)$ and $\ous(\mtn)$. 

We deduce that the $\osu$-algebra $\ous(\mtn)$ is generated by the
Chern classes $c_1(\rho(F))$, and is specified multiplicatively by
the ideal of relations amongst them. To compute this ideal, we
recall Davis and Januszkiewicz's space ${\it BP}^n$, which depends
only on the polytope $P^n$; all its cells are in even dimensions, and
the description of its cohomology extends immediately to complex
cobordism. Thus $\ous({\it BP}^n)$ is isomorphic to the Stanley-Reisner
$\osu$-algebra of $P^n$, which is the quotient of the polynomial
algebra $\osu[x_F:F\in\mathcal{F}]$ by a certain ideal $I$,
generated by those squarefree monomials
$\prod_\mathcal{E}x_F=x_\mathcal{E}$ for which $\cap_\mathcal{E}F$ is
empty. 

Since ${\it BP}^n$ is homotopy equivalent to the Borel construction
$ET^n\times_{T^n}\mtn$, there is a fibration
\begin{equation}\label{torfib}
T^n\rightarrow\mtn\xrightarrow{j}{\it BP}^n,
\end{equation}
classified by a map $l\colon {\it BP}^n\rightarrow BT^n$. The map $j$
pulls each cobordism class $x_F$ back to the Chern class
$c_1(\rho(F))$, whilst $l$ pulls the $i$th Chern class $c_1$ back to
some element $\lambda_i$, both in cohomology and complex cobordism,
for $1\leq i\leq n$. Considered as a homomorphism on $2$-dimensional
generators, we may identify
$l^*\colon\bZ^n\rightarrow\bZ^\mathcal{F}$ with the dual of the
dicharacteristic of $(\mtn,\alpha,\pi,P^n)$; in this setting we
abuse notation by interpreting $\lambda_i$ as an element of the
polynomial algebra $\ous[x_F:F\in\mathcal{F}]$.

We then compare the Serre spectral sequence
\[
H^*({\it BP}^n;\ous(T^n))\Longrightarrow\ous(\mtn)
\]
of \eqref{torfib} with the corresponding spectral sequence for the
universal principal fibration, which pulls back to the former along
$l$. The only differential is $d_2$, which annihilates all products of
one dimensional elements. We deduce the following result, to be
compared with the Danilov-Jurkiewicz Theorem \cite{daja:cpc}
describing the integral cohomology of toric varieties. 
\begin{prop}\label{cobstr}
Given any omnioriented toric manifold $(\mtn,\alpha,\pi,P^n)$, the
cobordism ring $\ous(\mtn)$ is isomorphic to
\[
\ous[x_F:F\in\mathcal{F}]\;/\;(I+J),
\]
where $J$ denotes the homogeneous ideal generated by the
$\lambda_i$, and the elements $x_F$ depend on the omniorientation.
\end{prop}

In Proposition \ref{cobstr}, each $x_F$ corresponds to the Chern class
$c_1(\rho(F))$. Reversing the facial orientation of a single $\xft$
therefore applies the inverse of the universal formal group law to
$x_F$. This is linked by Poincar\'e duality to the effect on bordism
theory, which manifests itself in a change of stably complex structure
on $\mtn$, and on those submanifolds $\xgt$ for which $F$ meets
$G$. The manifolds themselves remain unaltered.

In the case of bounded flag manifolds, we obtained a result equivalent
to Proposition \ref{cobstr} in \cite{bura:fml}. We did not, however,
specify an omniorientation there, but worked instead with the stably
complex structure of Example \ref{bnex}.

\section{Connected sums}\label{cs}

In order to construct connected sums of omnioriented compact toric
manifolds, we introduce an operation of connected sum for simple
polytopes equipped with extra combinatorial data. We work in dimensions
$\geq 2$, and deal separately with the degenerate case $n=1$ at the
end. Wherever practicable we write $m(P)$ for the number of facets of
$P^n$ and $q(P)$ for the number of vertices.

Before we begin we introduce a polyhedral template $\varGamma^n$,
which is the intersection of $n$ halfspaces in $\bR^n$. Strictly
speaking, it fails to qualify as a simple $n$-polyhedron because
$n<n+1$, but no contradiction arises from retaining the associated
terminology, and we do so for convenience. We embed the standard
$(n-1)$-simplex $\splx{n-1}$ in the subspace $\{x:x_1=0\}$ of
$\bR^{n-1}$, and construct $\varGamma^n$ by taking cartesian products
with the first coordinate axis. Its facets $G_r$ therefore have the form
$\bR\times D_r$, for $1\leq r\leq n$. Both $\varGamma^n$ and $G_r$ are
divided into positive and negative halves, determined by the sign of the
coordinate $x_1$.  

Given simple polytopes $P^n$ and $Q^n$ in $\bR^n$, we assume that
respective vertices $v$ and $w$ are distinguished. In addition, we order
the facets of $P^n$ meeting in $v$ as $E_r$, and the facets of $Q^n$
meeting in $w$ as $F_r$, for $1\leq r\leq n$. Recalling the notation of
Section \ref{scs}, we write $\mathcal{C}_v$ and $\mathcal{C}_w$ for the
complementary sets of facets; those in $\mathcal{C}_v$ avoid $v$, and
those in $\mathcal{C}_w$ avoid $w$. Their cardinalities are $m(P)-n$ and
$m(Q)-n$ respectively.

We now select a projective transformation $\phi_P$ which maps $v$ to
$x_1=+\infty$, and embeds $P^n$ in $\varGamma^n$ so as to satisfy two
conditions; firstly, that the hyperplane defining $E_r$ is identified
with the hyperplane defining $G_r$, for each $1\leq r\leq n$, and
secondly, that the images of the hyperplanes defining $\mathcal{C}_v$
meet $\varGamma^n$ in its negative half. This may be achieved, for
example, by considering the composition $T\cdot\phi_P'$, where $\phi'_P$
is an affine equivalence mapping $v$ and its vertex figure to
$(1,0,\dots,0)$ and $\splx{n}$ respectively, and $T$ is defined by
$T(x)=x/(1-x_1)$. We choose $\phi_Q$ similarly; it maps $w$ to
$x_1=-\infty$, and identifies the hyperplanes defining $F_r$ and $G_r$
in such a way that the images of the hyperplanes defining
$\mathcal{C}_w$ meet $\varGamma^n$ in its {\it positive\/} half. We
define the {\it connected sum} $P^n\cs_{v,w} Q^n$ of $P^n$ at $v$ and
$Q^n$ at $w$ to be the simple convex $n$-polytope determined by all
these hyperplanes. It is defined only up to combinatorial equivalence;
moreover, different choices for either of $v$ and $w$, or either of the
orderings for $E_r$ and $F_r$, are likely to affect the combinatorial
type. When the choices are clear, or their effect on the result
irrelevant, we use the abbreviation $P^n\cs Q^n$.

The face lattice $\mathfrak{L}_F(P\cs Q)$ is obtained from
$\mathfrak{L}_F(P)\cup\mathfrak{L}_F(Q)$ by identifying $E_r$ with $F_r$
for $1\leq r\leq n$; we write the result as $G_r$, and partition
the facets as
\begin{equation}\label{csfacets}
\mathcal{F}(P\cs Q)=
\mathcal{C}_v\cup\{G_r:1\leq r\leq n\}\cup\mathcal{C}_w.
\end{equation}
The vertices of $P^n\cs Q^n$ are the union of those of $P^n$ and $Q^n$,
omitting $v$ and $w$. Thus $m(P\cs Q)=m(P)+m(Q)-n$, and $q(P\cs
Q)=q(P)+q(Q)-2$.    

By way of illustration, we consider the connected sum $\splx{n}\cs_{v,w}
Q^n$, noting that the symmetry of the simplex guarantees that the result
is independent of the choice of $v$. We take $v$ to be $0$, so that
$\mathcal{C}_v=\{D_{n+1}\}$, and assume that $\phi_{\varDelta^n}$
identifies $D_r$ with $G_r$, for each $1\leq r\leq n$. So
$\phi_{\varDelta^n}(\varDelta^n)$ consists of $\varGamma^n$, truncated
in its negative half by a single hyperplane $H$ corresponding to the
image of $D_{n+1}$. Applying $\phi_Q^{-1}$, we deduce that the connected
sum is combinatorially equivalent to the polytope obtained from $Q^n$ by
including an extra hyperplane in the defining set. Such an $H$ must
isolate $w$, but no other vertex, and we interpret its inclusion as {\it
pruning} $Q^n$ at $w$. We write
\begin{equation}\label{pruopcs}
\splx{n}\cs_{v,w} Q^n\equiv\varPi_w(Q^n),
\end{equation}
where $\varPi_w$ denotes the appropriate pruning operator.

In order to generalise this example to products of simplices, we
need a pruning operator $\varPi_F$ for each face $F$ of $Q^n$; it is
defined in the obvious fashion, and detaches $F$ from $Q^n$ by any
hyperplane which separates the vertices of $F$ from the complementary
vertices of $Q^n$. Such operators obey two simple rules, which
lead to Theorem \ref{prucalc} below. Firstly, for any product
$P^m\times Q^n$ and any face $E$ of $P^m$, there are combinatorial
equivalences 
\begin{equation}\label{pruprod}
\begin{split}
\varPi_E(P^m)\times Q^n\equiv\varPi_{E\times Q}(P^m\times Q^n)
&\sands\\
&P^m\times\varPi_F(Q^n)\equiv\varPi_{P\times F}(P^m\times Q^n).
\end{split}
\end{equation}
Secondly, for any product of simplices
$\splx{m}\times\splx{n-m}$ with distinguished vertex $v$ there is a
face $G$ of $\splx{n}$, a vertex $v'$ not in $G$, and a combinatorial
equivalence 
\begin{equation}\label{prusplx}
\splx{m}\times\splx{n-m}\equiv\varPi_G(\splx{n})
\end{equation}
which maps $v$ to $v'$.

We now turn to arbitrary  products
$\splx{m_1}\times\dots\times\splx{m_k}$, where $m_1+\dots+m_k=n$.
\begin{thm}\label{prucalc}
Given any simple polytope $Q^n$, there is a combinatorial equivalence
\[
(\splx{m_1}\times\dots\times\splx{m_k})\cs_{v,w}Q^n\equiv
\varPi_{F_1}(\dots(\varPi_{F_k}(Q^n))\dots)
\]
for some sequence $F_i$ of products of simplices.
\end{thm}
\begin{proof}
We first extend \eqref{prusplx} to $k$-fold products by induction. For
$k\geq 3$, the inductive hypothesis provides a combinatorial equivalence
\[
\left(\splx{m_1}\times\dots\times\splx{m_{k-1}}\right)\times\splx{m_k}
\equiv
\varPi_{G_1}(\dots(\varPi_{G_{k-2}}(\splx{n-m_k}))\dots)\times\splx{m_k},
\]
where the $G_i$ are products of simplices; iterating \eqref{pruprod}
replaces the right-hand expression by
\[
\varPi_{G_1\times\varDelta^{m_k}}(\dots(\varPi_{G_{k-2}
\times\varDelta^{m_k}}(\splx{n-m_k}\times\splx{m_k}))\dots),
\]
and applying \eqref{prusplx} yields $\varPi_{F_1}(\dots
(\varPi_{F_{k-1}}(\splx{n})\dots)$, where the $F_i$ are products of
simplices as required. As in \eqref{prusplx}, we may ensure that $v$
corresponds to one of the original vertices of $\splx{n}$ (for which
we retain the label $v$) under the resulting equivalence.

Our connected sum is therefore combinatorially equivalent to
\begin{equation}\label{nearly}
\varPi_{F_1}(\dots(\varPi_{F_{k-1}}(\splx{n}))\dots)\cs_{v,w}Q^n.
\end{equation}
Since none of the faces $F_i$ contains the vertex $v$, they may be
identified with the corresponding faces of $\splx{n}\cs_{v,w}Q^n$; in
other words, \eqref{nearly} is combinatorially equivalent to
$\varPi_{F_1}(\dots(\varPi_{F_{k-1}}(\splx{n}\cs_{v,w}Q^n))\dots)$,
and the result follows by final appeal to \eqref{pruopcs}.
\end{proof}

We may now construct the connected sum of omnioriented toric
manifolds, assumed henceforth to be compact. Given
$(\mtn,\alpha,\pi,P^n)$ with fixed point $x$ projecting to the vertex
$v$ of $P^n$, and $(N^{2n},\beta,\mu,Q^n)$, with fixed point $y$
projecting to the vertex $w$ of $Q^n$, we suppose that the associated 
dicharacteristics are $\ell_M$ and $\ell_N$ respectively. We partition
the facets of $P^n\cs Q^n$ as in \eqref{csfacets}.
\begin{lem}\label{newdichs}
Up to $\theta$-translation, we may assume that $\ell_M$ identifies
$T^{F_r}$ with the $r$th coordinate subtorus $T_r$, for each $1\leq
r\leq n$.
\end{lem}
\begin{proof}
Since the subtori $T(F_r)$ generate $T^n$, we may define an
automorphism $\psi$ of $T^n$ by mapping $T(F_r)$ onto $T_r$,
preserving orientation, for each $1\leq r\leq n$. We conclude by
replacing $\ell_M$ with $\psi\cdot\ell_M$, and appealing to
Proposition \ref{oocop2}.
\end{proof}
Applying Lemma \ref{newdichs} to both $\ell_M$ and $\ell_N$ allows us to
combine them into a dicharacteristic 
\begin{equation}\label{ellcs}
\ell_\#(F)=
\begin{cases}
\ell_M(F)&\text{for $F\in\mathcal{C}_v$}\\
T_k&\text{for $F=G_r$ and $1\leq r\leq n$}\\
\ell_N(F)&\text{for $F\in\mathcal{C}_w$}
\end{cases}
\end{equation}
on $P^n\cs_{v,w}Q^n$.

We then define the equivariant connected sum
\[
(\mtn\cs_{x,y}N^{2n},\alpha\cs\beta,\pi\cs\mu,P^n\cs_{v,w}Q^n)
\]
to be the omnioriented toric manifold derived from \eqref{ellcs}.
Since $P^n\cs_{v,w}Q^n$ is determined only up to combinatorial
equivalence, we need Corollary \ref{uniquecs} to ensure that our
connected sum is well-defined. Its equivariant diffeomorphism type
depends on the choice of fixed points $x$ and $y$, as well as the
orderings of the facets $E_r$ and $F_r$. Nevertheless, Corollary
\ref{final} below shows that the properties we require of the induced
stably complex structure are suitably invariant.

\begin{thm}\label{cscs}
The manifold $\mtn\cs_{x,y}N^{2n}$ is diffeomorphic to the connected
sum of $\mtn$ and $N^{2n}$; furthermore, the diffeomorphism identifies
the stably complex structure induced by $\ell_\#$ with the connected
sum of those induced on $\mtn$ by $\ell_M$ and on $N^{2n}$ by
$\ell_N$.
\end{thm}
\begin{proof}
We consider the halfspace $H_\epsilon$ of $\bR^n$ given by $x_1\geq
\epsilon$ for some small $\epsilon>0$, and note that its inverse image
under the projective transformation $\phi_P$ is a halfspace $H_v$,
which intersects $P^n$ in a closed neighbourhood $N_v$ of the vertex
$v$. We then apply Lemma \ref{techlem} by choosing $P_1^n=P^n$,
$\;H_1=H_v$, $\;P_2^n=P^n\cs_{v,w}Q^n$, and $H_2=H_\#$; we take $f$ to
be $\phi_P\colon P_1^n\setminus N(v)\rightarrow 
P^n\cs_{v,w}Q^n\cap\{x_1<\epsilon\}$. We obtain an equivariant
diffeomorphism $\phi_P^+$ of $\mtn\setminus V_x$ into
$\mtn\cs_{x,y}N^{2n}$ (for some invariant neighbourhood $V_x$ of $x$),
which identifies the restrictions of the stably complex structures
induced by $\ell_M$ and $\ell_\#$ respectively. We then repeat the
process for $Q^n$, using the halfspace $x_1\leq-\epsilon$, and obtain
a corresponding diffeomorphism $\phi_Q^+$. The images of $\phi_P^+$
and $\phi_Q^+$ overlap in a collared $(2n-1)$-sphere, and the proof is
complete. 
\end{proof}
\begin{cor}\label{final}
For any choice of $x$, $y$, or orderings of the $E_r$ and $F_r$, the
stably complex structure induced on $\mtn\cs_{x,y} N^{2n}$ is cobordant
to the disjoint union of the structures induced on $\mtn$ and $N^{2n}$
respectively. 
\end{cor}
\begin{proof}
This is a direct consequence of Theorem \ref{cscs}, since the
connected sum of any two stably complex structures is cobordant to their 
disjoint union.
\end{proof}

Our main result follows from  Corollary \ref{final} and Example 
\ref{bijext}.
\begin{thm}\label{ccccc}
In dimensions $>2$, every complex cobordism class contains a toric
manifold, necessarily connected, whose stably complex structure is
induced by an omniorientation, and is therefore compatible with the
action of the torus.
\end{thm}
It follows from \eqref{dup} and Theorem \ref{prucalc} that the base
polyhedron of any such manifold may be constructed by pruning
$\splx{n}$ at an appropriate sequence of products of simplices, up to
combinatorial equivalence.

Finally, we address the degenerate case $n=1$, corresponding to toric
manifolds of dimension 2. Any simple polytope of dimension $1$ is a
$1$-simplex, and our definition of connected sum remains valid,
yielding $\varDelta_1^1\cs_{v,w}\varDelta_2^1=\splx{1}$. We cannot,
however, form the connected sum \eqref{ellcs} of dicharacteristics,
because the facets of $\splx{1}$ are its vertices, and the information
contained in $\ell_1(v)$ and $\ell_2(w)$ is lost when $v$ and $w$ are
deleted. We outline two different approaches to this anomaly, and
leave readers to decide on their preference.

We recall that $\varOmega^U_2$ is isomorphic to $\bZ$, and that the
cobordism class corresponding to a nonzero integer $n$ may be
represented by $\CPo$ with stably complex structure
$\tau\oplus\bR^2\cong\zeta(1)^n\oplus\zeta(1)^n$. According to Example
\ref{cpnex}, this structure is induced by an omniorientation only if
$n=\pm 1$. Nevertheless, $\CPo$ is a toric manifold, and the action of
the torus remains compatible with the exotic structure given by other
values of $n$. One option is therefore to assert Theorem \ref{ccccc}
in this weaker sense when $n=1$. Alternatively, we might retain the
values $\ell_1(v)$ and $\ell_2(w)$ on the connected sum by assigning
the information to the $1$-simplex (as opposed to its vertices). A
second option is therefore to extend the notion of dicharacteristic
when $n=1$, so that Theorem \ref{ccccc} holds as stated.

\end{document}